
  \magnification\magstep1
  \baselineskip = 0.5 true cm
                           
  \def\sa{\vskip 0.30 true cm}
  \def\sb{\vskip 0.60 true cm}
  \def\sc{\vskip 0.15 true cm}
  \def\sd{\vskip 0.50 true cm}
  \def\demi{ { {\lower 3 pt\hbox{$\scriptstyle 1$}} \over 
               {\raise 3 pt\hbox{$\scriptstyle 2$}} } } 

  \nopagenumbers
  \vsize = 25.7   true cm
  \hsize = 15.6   true cm
  \voffset = -1 true cm
  \parindent=0.85 true cm

  \rightline{December 1996} 

\sc
\sb
\sb 
\vskip 1 true cm  

\centerline {\bf ON QUADRATIC AND NONQUADRATIC FORMS: APPLICATION}

\vskip 0.8 true cm

\centerline {{\bf TO NONBIJECTIVE 
${\bf R}^{2m} \to {\bf R}^{2m-n}$ TRANSFORMATIONS}{\footnote{$^*$} 
{\sevenrm 
Work presented both to the Symposium ``Symmetries in Science IX'' 
held at the Cloister Mehrerau (Bregenz, Austria, 6-10 August 1996) 
and to the Symposium ``III Catalan Days on Applied Mathematics'' 
held at the Institut d'Estudis Ilerdencs (Lleida, Spain, 27-29 November 1996). 
To be published in {\bf Symmetries in Science IX}, 
ed.~Bruno Gruber (Plenum Press, New York, 1997).}}}

\vskip 0.4 true cm                                       

\sa
\sb

\vskip 0.5 true cm

\centerline 
{Maurice Kibler}
\sa

\centerline {Institut de Physique Nucl\'eaire de Lyon}
\centerline {IN2P3-CNRS et Universit\'e Claude Bernard}
\centerline {43 Boulevard du 11 Novembre 1918}
\centerline {F-69622 Villeurbanne Cedex}
\centerline {France}

\sa

\sa
\sa
\sb
\sb

\centerline {\bf Abstract} 

\sa
\sa

Hurwitz transformations are defined as specific automorphisms 
           of a Cayley-Dickson algebra. These transformations generate 
           quadratic and nonquadratic forms. We investigate here the 
           Hurwitz transformations corresponding to Cayley-Dickson algebras 
           of dimensions $2m = 2$, 4 and 8. The Hurwitz transformations 
           which lead to quadratic forms are discussed from geometrical and
           Lie-algebraic points of view. Applications to number theory and 
           dynamical systems are briefly examined.

\baselineskip = 0.7 true cm

  \sa

\vfill\eject

\baselineskip = 0.5 true cm

  \vglue 4 true cm

\noindent {\bf ON QUADRATIC AND NONQUADRATIC FORMS: APPLICATION}

\sb

\noindent {\bf TO NONBIJECTIVE 
${\bf R}^{2m} \to {\bf R}^{2m-n}$ TRANSFORMATIONS} 

\vskip 0.4 true cm 

\sa
\sb

\vskip 0.5 true cm

\leftskip = 1.6 true cm

{Maurice Kibler} 

\sa 

{Institut de Physique Nucl\'eaire de Lyon}

{IN2P3-CNRS et Universit\'e Claude Bernard}

{43 Boulevard du 11 Novembre 1918}

{F-69622 Villeurbanne Cedex}

{France}

\leftskip = 0 true cm

\sa
\sb
\sa
\sa
\sb
\sb

\baselineskip = 0.5 true cm

\noindent {\bf INTRODUCTION}

\sb

The application of (Hopf) fiber bundles is well developed 
in theoretical and mathematical physics.$^{1,2}$ 
Along this vein, the Hopf fibrations on spheres lead to
nonbijective canonical quadratic transformations useful in classical and
quantum mechanics.$^{3-13}$ The Hopf fibrations on spheres, as well as their
extensions on hyperboloids,$^{10}$ yield 
the concepts of `constraint Lie algebra'
and `Lie algebra under constraints' which are of importance for
connecting invariance or noninvariance algebras of dynamical systems.$^{11}$  

On the other hand, there is a close connection between number theory and
quadratic mappings, especially those mappings sending the ($2m - 1$)-sphere on 
the ($2m-1-n$)-sphere.$^{14,15}$ In this direction, one can easily imagine that
the replacement of mappings on spheres by mappings on hyperboloids might
produce interesting results. 

It is the aim of this work to describe how fibrations on spheres and
hyperboloids naturally arise in the framework of Cayley-Dickson algebras. 
This will provide us with quadratic$^{10}$ and nonquadratic$^{12}$ forms which 
generate nonbijective transformations from ${\bf R}^{2m}$ onto 
${\bf R}^{2m - n}$. We shall give two prototypical applications of these
quadratic and nonquadratic transformations, viz. one to number theory (in the
spirit of Ref.~15) and one to dynamical systems. Much of this work takes its
origin in Refs.~5 to 13. The application to number theory constitutes the
premises of a more complete work. 

  
\vfill\eject

\noindent {\bf CAYLEY-DICKSON ALGEBRAS} 

\sb

The Cayley-Dickson (C-D) algebras generalize the algebras of
complex numbers, quaternions and octonions. A $2m$-dimensional C-D algebra may 
be obtained from a $m$-dimensional ($m=1$, 2, 4, $\cdots$) C-D algebra by a 
`doubling' process (cf. ${\bf  C}={\bf  R} + {\rm i}{\bf  R}$). 
In the present paper, we limit ourselves to C-D algebras of dimensions 
$2m=2$, 4 and 8. They are denoted as $A(c)$ with 
$c \equiv c_1$, 
$c \equiv c_1,c_2$ and 
$c \equiv c_1,c_2,c_3$ 
for $2m=2$, $4$ and 8, respectively, where 
$c_i = \pm 1$ for $i=1,2,3$. Note that $A(-1      )$, $A(-1,-1   )$ and 
$A(-1,-1,-1)$ are nothing but the algebra ${\bf C}$ of usual complex numbers, 
the algebra ${\bf H}$ of ordinary (or Hamilton) quaternions and the algebra 
${\bf O}$ of ordinary (or Cayley) octonions, 
respectively. The algebra $A(1)$ is the algebra ${\bf \Omega}$ 
of hyperbolic complex numbers. The three algebras 
$A(c_1, c_2)$ with $(c_1,c_2) \not = (-1,-1)$ are isomorphic to the algebra 
${\bf  N}_1$ of hyperbolic quaternions and the seven algebras
$A(c_1,c_2,c_3)$ with $(c_1,c_2,c_3) \not = (-1,-1,-1)$ are isomorphic to the
algebra ${\bf  O}'$ of hyperbolic octonions. 

Let 
$$
u = u_0 + \sum_{k=1}^{2m - 1} u_k \> e_k
\eqno (1)
$$
be an element of $A(c)$ where the real numbers 
$u_0,u_1,\cdots,u_{2m-1}$ are the components of 
$u$ and the set 
$\{e_1,e_2,\cdots,e_{2m-1}\}$ is a system of generators of $A(c)$. The
product $w=uv$ of the hypercomplex numbers $u$ and $v$ of $A(c)$ is defined 
once the multiplication rule
$$
e_ie_j = - g_{ij} + \sum_{k=1}^{2m-1} a^k_{ij} \> e_k
\qquad i \ {\rm and} \ j = 1,2,\cdots,2m-1
\eqno (2)
$$
for the generators of $A(c)$ is given. The constants $a^k_{ij}$ 
(totally antisymmetric in $ijk$) appear in Ref.~10 for $2m=2$, 
$4$ and 8. The constants $g_{ij}$ are defined by the matrix
$$
g \equiv (g_{ij}) 
= {\rm diag} (-c_1, -c_2, c_1c_2, -c_3, c_1c_3, c_2c_3, -c_1c_2c_3)
\eqno (3)
$$
for $2m=8$. The matrix $g$ for $2m=4$ (or $2m =2$) follows by restricting 
the matrix in Eq.~(3) to its first four (or two) lines and columns.

The product $w=uv$ can be represented in a matrix form as
$$
{\bf {w}} = H({\bf {u}} ; c) \> {\bf {v}}
\eqno (4)
$$
where $\bf {w}$, $\bf {u}$ and $\bf {v}$ are column vectors defined by
$$
^t{\bf {w}} = (w_0,w_1,\cdots,w_{2m-1}) \ 
^t{\bf {u}} = (u_0,u_1,\cdots,u_{2m-1}) \
^t{\bf {v}} = (v_0,v_1,\cdots,v_{2m-1})
\eqno (5)
$$ 
and $H({\bf {u}} ; c)$ is a $2m \times 2m$ matrix generalizing the Hurwitz
matrix$^{\rm 10}$ that occurs in the search of bilinerar forms 
$w_\alpha = \phi_\alpha (u_\beta,v_\beta)$ with $\alpha = 1,2,\cdots,n$ 
                                            and $\beta  = 1,2,\cdots,n$ such
that
$$
 w^2_1 + w^2_2 + \cdots + w^2_n =
(u^2_1 + u^2_2 + \cdots + u^2_n) \> 
(v^2_1 + v^2_2 + \cdots + v^2_n) 
\eqno (6)
$$
The matrix $H({\bf {u}} ; c)$ depends on $c$ and the components 
of ${\bf u}$. For $2m=8$, the matrix $H({\bf {u}} ; c)$ is the 
following one 
$$
\pmatrix {
u_0 & c_1u_1 & c_2u_2 & -c_1c_2u_3 & c_3u_4 & -c_1c_3u_5 &-c_2c_3u_6
&c_1c_2c_3u_7\cr
\cr
u_1 & u_0 & c_2u_3 & - c_2u_2 & c_3u_5 & -c_3u_4 & c_2c_3u_7 & -c_2c_3u_6\cr
\cr
u_2 & -c_1u_3 & u_0 &  c_1u_1 & c_3u_6 & -c_1c_3u_7 & -c_3u_4 & c_1c_3u_5\cr
\cr
u_3 & -u_2 & u_1 & u_0 & c_3u_7 & -c_3u_6 & c_3u_5 & -c_3u_4\cr
\cr
u_4 & -c_1u_5 & -c_2u_6 & c_1c_2u_7 & u_0 & c_1u_1 & c_2u_2 & -c_1c_2u_3\cr
\cr
u_5 & -u_4 & -c_2u_7 & c_2u_6 & u_1 & u_0 & -c_2u_3 & c_2u_2\cr
\cr
u_6 & c_1u_7 & -u_4 & -c_1u_5 & u_2 & c_1u_3 & u_0 & -c_1u_1\cr
\cr
u_7 & u_6 & -u_5 & u_4 & u_3 & u_2 & -u_1 &u_0\cr
}
\eqno (7)
$$
The matrix $H({\bf {u}} ; c)$ for $2m=4$ can be obtained from 
(7) by omitting the 5-, 6-, 7- and 8-th lines and columns. 
Similarly, the matrix $H({\bf {u}} ; c)$
for $2m=2$ correponds to the two first lines and two first columns 
of the matrix (7). It can be shown that the matrix $H({\bf {u}} ; c)$ 
for the algebra $A(c)$ 
can be developed as a linear combination of Clifford matrices. Indeed, we have 
$$
H({\bf u} ; c) = u_0 {\bf 1} + \sum_{k=1}^{2m-1} u_k \; {^t\Gamma_k}
\eqno (8)
$$
where {\bf 1} is the $2m \times 2m$ unit matrix and the $2m \times 2m$
(Clifford) matrices $\Gamma_k$ satisfy 
$$
\Gamma_i \Gamma_j + \Gamma_j \Gamma_i = - 2 \> g_{ij} \> {\bf 1}
\qquad (i \; \hbox { and } \; j = 1, 2, \cdots, 2m-1)
\eqno (9)
$$
It is thus
possible to associate a Clifford algebra ${\cal C}(p,q)$ of 
degree $p+q = 2m-1$ (with
$q - p$ being the signature of $g$) to the $2m$-dimensional C-D algebra $A(c)$.
An important property of the matrix $H({\bf {u}} ; c)$ is the following.

\sb

{\bf Property 1}. The matrix $H({\bf {u}} ; c)$ belongs to 
${\bf  R}^+ \times {\rm O}(2m)$ or 
${\bf  C}   \times {\rm O}(m,m)$ according to as the $2m \times 2m$ metric
$$
\eta = 1 \oplus g
\eqno (10)
$$
has the signature $2m$ of 0. More specifically, we have 
$$
^tH({\bf {u}} ; c) \> \eta \> 
  H({\bf {u}} ; c) = (^t{\bf {u}} \> \eta \> {\bf {u}}) 
\> \eta
\eqno (11)
$$

We shall use the qualification compact or noncompact according to as the metric
$\eta$ is Euclidean or pseudoEuclidean. The C-D algebra $A(c)$ may be normed in
the compact case or pseudonormed (corresponding to a singular division algebra
with a cone of zero divisors) in the noncompact case. 
We close this introductory section with a lemma that is of central
importance for what follows.

\sb

{\bf Lemma}. The complex conjugation  
in   a $2m$-dimensional C-D algebra $A(c)$ and
in its  $m$-dimensional C-D subalgebras induces the existence of 
$2m - \delta(m,1)$ anti-involutions on $A(c)$. 

The latter anti-involutions (i.e. involutive anti-automorphisms) 
$j : A(c) \rightarrow A(c)$ can be described 
in the following way. In the case $m \not = 1$,
one anti-involution $j_0$ corresponds to the complex conjugation in the C-D
algebra $A(c)$ while the $2m-1$ remaining anti-involutions
$j_1,j_2,\cdots,j_{2m-1}$ correspond to the complex conjugation in the various 
$m$-dimensional C-D subalgebras of $A(c)$. (The complex conjugate $\overline u$
of the element $u$ of Eq.~(1) is defined in $A(c)$ by
     $\overline u=u_0 - \sum_{k=1}^{2m-1} u_ke_k$.) For $m=1$, the 
two-dimensional C-D algebras 
$A(c_1)$ have only one anti-involution since $j_0$ coincides with $j_1$. 

\sb 
\sd

\noindent {\bf HURWITZ TRANSFORMATIONS}
                    
\sb

We associate to the element $u$ of $A(c)$, see Eq.~(1), 
the hypercomplex number
$$
\widehat {u} = u_0 + \sum_{k=1}^{2m-1} \varepsilon_k \> u_k \> e_k
\eqno (12)
$$
where $\varepsilon_k = \pm 1$ for $k = 1,2,\cdots,2m-1$. The notation
$$
{\bf \varepsilon} = {\rm diag} (1,\varepsilon_1,\cdots,\varepsilon_{2m-1})
\eqno (13)
$$
shall be used below.

\sb

{\bf Definition 1}. The mapping
$$
T[1;c;{\bf \varepsilon}] : A(c) \rightarrow A(c) 
                         : u \mapsto x = u \> \widehat {u}
\eqno (14)
$$
is called a (right) Hurwitz transformation. 

The Hurwitz transformations can be classified into several types 
according to the various
possibilities for the $2m \times 2m$ matrix ${\bf \varepsilon}$.

{\it Type} $A_1$. For $\widehat u = u$: The transformations 
$T[1;c;{\bf \varepsilon}]$, where ${\bf \varepsilon}$ is the unit matrix 
{\bf 1}, were called quasiHurwitz transformations in Ref.~10 and 
transformations of type $A_1$ in Ref.~12. 

{\it Type} $B_1$. For ${\widehat u}=j(u)$: The transformations 
$T[1;c;{\bf \varepsilon}]$, 
where ${\bf \varepsilon}$ corresponds to an arbitrary anti-involution $j$ of 
$A(c)$, were called Hurwitz transformations in Ref.~10 and transformations 
of type $B_1$ in Ref.~12.

{\it Type} $C_1$. For $\widehat {u} \not = u$ or $j(u)$: 
The transformations $T[1;c;{\bf \varepsilon}]$, 
                    where ${\bf \varepsilon}$ corresponds neither to the 
unit matrix 
{\bf 1} nor to an anti-involution $j$ of $A(c)$, were called pseudoHurwitz
transformations in Ref.~10 and transformations of type $C_1$ in Ref.~12. 

\sb 
\sd

\noindent {\bf QUADRATIC TRANSFORMATIONS}

\sb

The mapping (14) can be rewritten in matrix form as
$$
{\bf R}^{2m \times 1} \rightarrow 
{\bf R}^{2m \times 1} : {\bf u} \mapsto {\bf x} = H ({\bf u} ; c) \> 
{\bf \varepsilon} \> {\bf u}
\eqno (15)
$$
The column vector {\bf x} contains $n$ components that are equal to $0$. The
integer $n$, which can take several values $(0 \leq n \leq 2m -1)$, depends on
${\bf \varepsilon}$. In other words, the Hurwitz transformation $T[1;c;{\bf
\varepsilon}]$ generates, via Eq.~(15), a ${\bf R}^{2m} \rightarrow {\bf
R}^{2m-n}$ mapping: To each vector, of components 
$u_\alpha$ ($\alpha = 0,1,\cdots,2m-1$), in ${\bf R}^{2m}$ is associated a 
vector in ${\bf R}^{2m-n}$
the components of which are quadratic functions of the variables $u_\alpha$.
Let us emphasize the following property which directly follows from Eq.~(11).

\sb

{\bf Property 2}. The relation
$$
      {^t {\bf x}} \> \eta \> {\bf x} = 
\big ({^t {\bf u}} \> \eta \> {\bf u} \big )^2
\eqno (16)
$$
is valid for all the Hurwitz transformations $T[1;c;{\bf \varepsilon}]$.

We continue with some examples for $2m = 8$, $4$ and 2. 
The ${\bf R}^{2m} \rightarrow {\bf R}^{2m-n}$ 
quadratic transformations shall be given in
explicit form only for $2m=8$. The cases $2m=4$ and $2m=2$ can be deduced from
the case $2m=8$ by taking $u_4=u_5=u_6=u_7=0$ with $c_3=0$ and
$u_2=u_3=u_4=u_5=u_6=u_7=0$ with $c_2=c_3=0$, respectively.

{\it Type} $A_1$. The mapping $T[1;c_1,c_2,c_3;{\bf 1}]$ gives the 
generic ${\bf R}^8 \rightarrow {\bf R}^8$ quadratic transformation
$$
x_0=u^2_0+c_1u^2_1+c_2u^2_2-c_1c_2u^2_3
       +c_3u^2_4-c_1c_3u^2_5-c_2c_3u^2_6+c_1c_2c_3u^2_7
$$ 
$$
x_1 =2u_0u_1 \qquad
x_2 =2u_0u_2 \qquad
x_3 =2u_0u_3 \qquad
x_4 =2u_0u_4 
$$
$$
x_5 =2u_0u_5 \qquad
x_6 =2u_0u_6 \qquad
x_7 =2u_0u_7 
\eqno (17)
$$
with the property
$$
  x^2_0 - c_1 x_1^2 - c_2 x_2^2 + c_1c_2 x_3^2 - c_3 x^2_4 + 
  c_1c_3 x^2_5 + c_2c_3 x^2_6 - c_1c_2c_3 x^2_7 = 
$$
$$
=(u_0^2 - c_1 u_1^2 - c_2 u_2^2 + c_1c_2 u_3^2 - c_3 u_4^2 + 
  c_1c_3 u^2_5 + c_2c_3 u^2_6 - c_1c_2c_3 u^2_7)^2 
\eqno (18)
$$ 
Equation (17) can be written as
$$
x_0 = 2u_0^2 - {^t{\bf u}} \> \eta \> {\bf u} \qquad 
x_k = 2u_0u_k \qquad (k = 1,2,\cdots,7)
\eqno (19)
$$
in condensed form.

From the geometrical point of view, the ${\bf R}^{2m} \rightarrow {\bf R}^{2m}$
transformations ($2m=8$, $4$ and 2) may be arranged in two classes. The compact 
case corresponds to the fibration 
$S^{2m-1} \rightarrow S^{2m-1}/Z_2 \sim {\bf R}P^{2m-1}$ of discrete fiber 
$Z_2$ while the noncompact case corresponds
to the fibration $H^{2m-1}(m,m) \rightarrow H^{2m-1}(m,m)/Z_2$ of
discrete fiber $Z_2$. (We use  $H^{2m-1}(m,m)$  to denote the hyperboloid of
equation $\sum_{i=1}^{m} u_i^2 - \sum_{i=m+1}^{2m} u_i^2$ in ${\bf R}^{2m}$.)

Note that the transformation $T[1;-1;{\bf 1}]$ corresponds to the so-called
Levi-Civita (conformal) transformation used in the restricted three-body 
problem.

{\it Type} $B_1$. The mapping 
$T[1;c_1,c_2,c_3;{\bf \varepsilon}]$ where
$$
{\bf \varepsilon} = {\rm diag} (1,-1,1,1,1,1,-1,-1)
\eqno (20)
$$
gives the generic ${\bf R}^8 \rightarrow {\bf R^5}$ quadratic transformation
$$
x_0 =u_0^2-c_1u_1^2+c_2u_2^2-c_1c_2u_3^2+c_3u_4^2-c_1c_3u_5^2
     +c_2c_3u_6^2-c_1c_2c_3u_7^2
$$
$$
x_2 =2\big ( u_0u_2+c_1u_1u_3+c_3u_4u_6-c_1c_3u_5u_7 \big ) 
$$
$$
x_3 =2\big ( u_0u_3+u_1u_2-c_3u_5u_6+c_3u_4u_7 \big )
$$
$$
x_4 =2\big ( u_0u_4+c_1u_1u_5-c_2u_2u_6+c_1c_2u_3u_7 \big ) 
$$
$$
x_5 =2\big ( u_0u_5+u_1u_4+c_2u_3u_6-c_2u_2u_7 \big ) 
\eqno (21)
$$
with the property
$$
x^2_0 - c_2 x_2^2 + c_1c_2 x_3^2 - c_3 x_4^2 + c_1c_3 x_5^2 = 
$$ 
$$
=(u_0^2 - c_1u_1^2 - c_2u_2^2 + c_1c_2u_3^2 - c_3u_4^2 + c_1c_3u^2_5
                              + c_2c_3u^2_6 - c_1c_2c_3u^2_7)^2 
\eqno (22) 
$$ 
The matrix ${\bf \varepsilon}$ in Eq.~(20) corresponds to an anti-involution 
of 
$A(c_1,c_2,c_3)$ 
of type $j_1,j_2,\cdots,j_7$. The other anti-involutions of type
$j_1,j_2,\cdots,j_7$ would lead to transformation formulae equivalent to 
Eqs.~(21) and (22). 

The mapping $T[1;c_1,c_2,c_3;{\bf \varepsilon}]$ where
$$
{\bf \varepsilon} = {\rm diag} (1,-1,-1,-1,-1,-1,-1,-1)
\eqno (23)
$$
gives the generic ${\bf R}^8 \rightarrow {\bf R}$ quadratic transformation
$$
x_0 
=u_0^2 - c_1u_1^2 - c_2u_2^2 + c_1c_2u_3^2 - 
         c_3u_4^2 + c_1c_3u^2_5 + c_2c_3u^2_6 - c_1c_2c_3u^2_7 
\eqno (24)
$$
This transformation (which reads $x_0 = {^t{\bf u}} \> \eta \> {\bf u}$ 
in a more
compact form) corresponds to the anti-involution $j_0$ of $A(c_1,c_2,c_3)$, 
i.e. to ${\widehat u} = j_0(u) = {\overline u}$. 

\sb

\leftskip = 0.7 true cm 
\rightskip= 0.7 true cm 

\noindent 
         {\bf Table 1.} 
         Fibrations (up to homeomorphisms) and Lie algebras under 
         constraints for the  
         ${\bf R}^{2m} \rightarrow {\bf R}^{2m-n}$ 
         quadratic transformations with $n=m-1+\delta (m,1)$.

\def\init{\tabskip 0pt\offinterlineskip}
\def\crr{\cr\noalign{\hrule}}

$$\vbox{\init\halign to 14.2 true cm{\strut#&\vrule#%
\tabskip=0em plus 2 em&
\hfil$#$\hfil&
\vrule#&
\hfil$#$\hfil&
\vrule#&
\hfil$#$\hfil&
\vrule#&
\hfil$#$\hfil&
\vrule#&
\hfil$#$\hfil&
\vrule#&
\hfil$#$\hfil&
\vrule#\tabskip 0pt\crr
&& && && && && && &\cr
&&\hbox {Transfor-} &&\hbox {Fibration} 
&&\hbox {Fiber} &&L &&L_0 &&L_1 &\cr
&&\hbox {mation} && && && && && &\cr
&& && && && && && &\crr
&& && && && && && &\cr
&& &&S^1 \rightarrow \{ 1 \} 
&&S^1 && &&{\rm so}(2) &&{\rm so}(2,1) &\cr
&&{\bf R}^2\rightarrow {\bf R} && && &&{\rm sp}(4,{\bf R}) && && &\cr
&& &&{\bf R} \rightarrow \{ 1 \} &&{\bf R} && &&{\rm so}(1,1) &&{\rm so}(2,1) &\cr
&& && && && && && &\crr
&& && && && && && &\cr
&& &&S^3 \rightarrow S^2 &&S^1 && &&{\rm so}(2) &&{\rm so}(4,2) &\cr
&&{\bf R}^4\rightarrow {\bf R}^3 
&&{\bf R}^2 \times S^1 \rightarrow {\bf R}^2 &&S^1 
&&{\rm sp}(8,{\bf R}) &&{\rm so}(2) &&{\rm so}(4,2) &\cr
&& &&{\bf R}^2 \times S^1 \rightarrow {\bf R} \times S^1 
&&{\bf R} && &&{\rm so}(1,1) &&{\rm so}(3,3) &\cr
&& && && && && && &\crr
&& && && && && && &\cr
&& &&S^7 \rightarrow S^4 &&S^3 && &&{\rm so}(3) &&{\rm so}(6,2) &\cr
&&{\bf R}^8\rightarrow {\bf R}^5  
&&{\bf R}^4 \times S^3\rightarrow {\bf R}^4 &&S^3 
&&{\rm sp}(16,{\bf R}) &&{\rm so}(3) &&{\rm so}(6,2) &\cr
&& &&{\bf R}^4 \times S^3 \rightarrow {\bf R}^2 \times S^2 
&&{\bf R}^2\times S^1 && &&{\rm so}(2,1) &&{\rm so}(4,4) &\cr
&& && && && && && &\crr
}}$$

\leftskip = 0 true cm 
\rightskip= 0 true cm 

\sa

The geometrical classification of the ${\bf R}^{2m} \rightarrow {\bf R}^{m+1}$
transformations ($2m=8$ and 4) associated to an anti-involution of type
$j_1,j_2,\cdots,j_{2m-1}$ of $A(c)$ leads to three classes. In the compact case,
we have the fibration on spheres $S^{2m-1} \rightarrow S^m$ of compact fiber
$S^{m-1}$. The noncompact case corresponds to two kinds of fibrations on
hyperboloids, viz. one with a compact fiber and the other one with a
noncompact fiber.$^{10}$ In addition, 
the ${\bf R}^{2m} \rightarrow {\bf R}$ transformations ($2m=8$, $4$ and 2)
associated to the anti-involution $j_0$ of $A(c)$ may be classified 
in two classes. The compact case corresponds to the fibration
$S^{2m-1} \rightarrow \{1\}$ of compact fiber $S^{2m-1}$ and the noncompact
case to the fibration $H^{2m-1}(m,m) \rightarrow \{ 1 \}$ of noncompact fiber
${\bf R}^m\times S^{m-1}$ (up to homeomorphisms). The 
fibrations corresponding to the anti-involution
$j_1$ ($\equiv j_0$) for $2m=2$ and to the anti-involutions of type
$j_1,j_2,\cdots,j_{2m-1}$ for $2m=4$ and 8 are reported in Table~1. 

The case $2m=4$ deserves three remarks. First, is to be noted that the
transformation $T[1;-1,-1;{\bf \varepsilon}]$ with 
${\bf \varepsilon} = {\rm diag} (1,-1,1,1)$ 
can be identified to the so-called Kustaanheimo-Stiefel$^{\rm 3}$ 
transformation (associated to the celebrated Hopf fibration on spheres 
$S^3 \rightarrow S^2 \sim S^3/S^1 \sim {\bf C}P^1)$ 
used for the regularization of the Kepler problem in classical mechanics.
Furthermore, the transformation $T[1;-1,1;{\bf \varepsilon}]$ with ${\bf
\varepsilon} = {\rm diag} (1,-1,1,1)$ 
corresponds to the transformation introduced by
Iwai$^{\rm 7}$ for reducing an Hamiltonian system by an $S^1$ action. 
Finally, the
transformation $T[1;c_1,c_2;{\bf \varepsilon}]$ with 
${\bf \varepsilon} = {\rm diag} (1,-1,1,1)$ 
and $c_1=-c_2=1$ (or $c_1=c_2=1$) corresponds to a transformation,
introduced by Lambert and Kibler,$^{\rm 10}$ inequivalent to the two previous
ones.

{\it Type} $C_1$. 
The mapping $T[1;c_1,c_2,c_3;{\bf \varepsilon}]$, where ${\bf \varepsilon}$ 
neither is the unit matrix ${\bf 1}$ nor corresponds to an anti-involution
$j$ of $A(c_1,c_2,c_3)$, yields new 
quadratic transformations only when $\sum_{k=1}^7 \varepsilon_k = -3$ or $5$. 
For example, the case 
$$
{\bf \varepsilon} = {\rm diag} (1,-1,-1,-1,-1,-1,1,1)
\eqno (25)
$$
gives the generic ${\bf R}^8 \rightarrow {\bf R}^7$ quadratic transformation
$$
x_0  = u_0^2 - c_1 u_1^2 - c_2 u_2^2 + c_1c_2 u_3^2
     - c_3 u_4^2 + c_1c_3 u_5^2 - c_2c_3 u_6^2 + c_1c_2c_3 u_7^2
$$
$$
x_2  = 2 (-c_3u_4u_6 + c_1c_3u_5u_7) \qquad 
x_3  = 2 (-c_3u_4u_7 +    c_3u_5u_6) 
$$
$$
x_4  = 2 ( c_2u_2u_6 - c_1c_2u_3u_7) \qquad
x_5  = 2 ( c_2u_2u_7 -    c_2u_3u_6) 
$$
$$
x_6  = 2 (    u_0u_6 -    c_1u_1u_7) \qquad 
x_7  = 2 (    u_0u_7 -       u_1u_6) 
\eqno (26)
$$
with the property
$$
x_0^2 - c_2 x_2^2 + c_1c_2 x_3^2 - c_3 x_4^2 + c_1c_3 x_5^2 + c_2c_3 x_6^2
- c_1c_2c_3 x_7^2 = 
$$
$$
=(u_0^2-c_1u_1^2-c_2u_2^2+c_1c_2u_3^2-c_3u_4^2+c_1c_3u^2_5+c_2c_3u^2_6
-c_1c_2c_3u^2_7)^2
\eqno (27)
$$
Other choices for ${\bf \varepsilon}$ with
$\sum_{k=1}^7 \varepsilon_k = -3$ or 5 lead to transformation formulae 
equivalent to Eqs.~(26) and (27). Remark that the restriction of 
Eq.~(26) to
the cases $2m=4$ and $2m=2$ does not lead to new transformations: The obtained
quadratic transformations are Hurwitz transformations of type $B_1$.
Furthermore, the (true) transformations of type $C_1$ are equivalent to
transformations of type $A_1$ for $2m=2$ and $2m=4$. 

For $2m=8$, the ${\bf R}^8 \rightarrow {\bf R}^7$ transformations 
provide explicit realizations for (i) the 
Hopf fibration $S^7 \rightarrow {\bf C}P^3$ of compact fiber $S^1$ when
$c_1+c_2+c_3 = -3$ and for (ii) its noncompact analogues, namely, 
${\bf R}^4 \times S^3 \rightarrow {\bf R}^4 \times S^2$ of compact fiber
${S}^1$ and 
${\bf R}^4 \times S^3 \rightarrow {\bf R}^3 \times S^3$ of noncompact fiber
${\bf R}$ when $c_1+c_2+c_3 \not = -3$. 


\vfill\eject

\noindent {\bf DIFFERENTIAL ASPECTS OF QUADRATIC TRANSFORMATIONS}

\sb

Differential and Lie-algebraic aspects of transformations of type $A_1$ and
$B_1$ were studied in Refs.~10 and 11. Here, we briefly discuss
the situation for the ${\bf R}^{2m} \rightarrow {\bf R}^{2m-n}$ quadratic
transformations of type $B_1$ with $n=m-1+\delta(m,1)$ and $2m=2$, $4$ and 8.

Let us consider the $2m$ one-forms defined by the column vector
$$
\Omega = 2 \> H({\bf u};c) \> {\bf \varepsilon} \> d{\bf u}
\eqno (28)
$$
where the $2m \times 2m$ matrix ${\bf \varepsilon}$ is associated to a given
anti-involution, of $A(c)$, of the type $j_1, j_2, \cdots, j_{2m-1}$.
Equation (28) provides us with (i) $2m-n$ total differentials that may be
integrated to give the nonvanishing components $x_j$ of ${\bf x}$ in Eq.~(15)
and (ii) $n$ one-forms $\omega_i$ that correspond to the vanishing components
of ${\bf x}$ in Eq.~(15). As a corollary of Property~1, we obtain
$$
^t\Omega \> \eta \> \Omega = 4 \; (^t {\bf u} \> \eta \>  {\bf u}) \;
                                  (^td{\bf u} \> \eta \> d{\bf u})
\eqno (29)
$$
In view of the nonbijective character of the 
${\bf R}^{2m} \rightarrow {\bf R}^{2m-n}$ transformation, we assume that each
one-form $\omega_i$ satisfies the constraint condition $\omega_i = 0$. The
introduction of these $n$ conditions in 
Eq.~(29) makes it possible to connect the
line element in ${\bf R}^{2m}$, with the metric $\eta$, and a line element in
                ${\bf R}^{2m-n}$.

To each one-form $\omega_i$, we may associate a vector field $X_i$, defined
in the symplectic Lie algebra ${\rm sp}(4m,{\bf R})$, with an action of 
$\omega_i$ on
${1 \over 2r} X_j$, where $r = {^t{\bf u}} \eta {\bf u}$, such that 
$\omega_i [{1 \over 2r} X_j] = \delta(i,j)$. Each vector field $X_i$ has
the property that $X_i \psi = 0$ for any function $\psi$ in 
$C^1({\bf R}^{2m-n})$. The following theorems and definitions are interesting 
for physical applications. 

\sb

{\bf Theorem 1}. The $n=m-1+\delta(m,1)$ vector fields $X_i$ span a Lie
algebra $L_0$ with respect to the commutator law.

\sb

{\bf Definition 2}. The subalgebra $L_0$ of $L = {\rm sp}(4m, {\bf R})$ 
is called a {\it constraint Lie algebra}. Let us define the subalgebra 
$L_1$ of $L$ by 
$$
L_1 = {\rm cent}_L L_0/L_0
\eqno (30)
$$
when $L_0$ is one-dimensional and by 
$$
L_1 = {\rm cent}_L L_0 
\eqno (31)
$$
when $L_0$ is semisimple. The Lie algebra $L_1$ is called a 
{\it Lie algebra under constraints}.

\sb

{\bf Theorem 2}. The constraint Lie algebra $L_0$ and the Lie algebra 
under constraints $L_1$ are characterized by the compact or noncompact nature 
of the fiber of the fibration associated to the corresponding 
${\bf R}^{2m} \rightarrow {\bf R}^{m+1-\delta(m,1)}$ quadratic transformation.

The triples $(L,L_0,L_1)$, $L_0 \subset L_1 \subset L$, 
are reported in Table~1 
for the quadratic transformations ${\bf R}^2 \rightarrow {\bf R}$, 
                                  ${\bf R}^4 \rightarrow {\bf R}^3$ and 
                                  ${\bf R}^8 \rightarrow {\bf R}^5$.


\vfill\eject

\noindent {\bf NONQUADRATIC TRANSFORMATIONS}

\sb

Nonquadratic transformations may be obtained by replacing Eq.~(15) by
$$
{\bf R}^{2m \times 1} \rightarrow 
{\bf R}^{2m \times 1} : 
{\bf u} \mapsto {\bf x} = H({\bf u};c)^{N} \> {\bf \varepsilon } \> {\bf u}
\quad \hbox {with} \quad N \in {\bf Z} - \{1\}
\eqno (32)
$$
Equations (15) and (32) define ${\bf R}^{2m} \rightarrow {\bf R}^{2m-n}$ 
quadratic and nonquadratic transformations ($0 \leq n \leq 2m-1$) in a 
unified way. The property (16) may be generalized as
$$
 ^t{\bf x} \> \eta \> {\bf x} = 
(^t{\bf u} \> \eta \> {\bf u})^{N+1}
\eqno (33)
$$
for the nonquadratic transformations. The quadratic and nonquadratic
transformations are associated to Hurwitz transformations that we shall denote
as $T[N;c;{\bf \varepsilon}]$ with $N \in {\bf Z}$. They can be
classified into transformations of type $A_N$, $B_N$ and $C_N$ by
employing the same character of distinction as for the transformations of type
$A_1$, $B_1$ and $C_1$. We limit here our consideration to two specific
examples: (1) transformations of type $A_N$ with $N \in {\bf Z} - \{-1\}$ and
$2m = 2$ and (2) transformations of type $B_{-1}$ with $2m = 4$.

\sb

{\bf Example 1}. The mapping $T[N;c_1;{\bf 1}]$ defines the
transformations ${\bf u} \mapsto {\bf x} = H({\bf u} ; c_1)^{N} \> {\bf u}$.
General properties of these transformations of type $A_N$ ($N \not = -1$) 
can be easily derived for $2m=2$. Indeed, we have the five following
properties: 
$$
 x_0^2 - c_1x_1^2 = 
(u_0^2 - c_1u_1^2)^{N+1}
$$
$$
d{\bf x} = (N+1) \> H({\bf u};c_1)^N \> d{\bf u}
$$
$$
dx_0^2 -c_1dx_1^2 = (N+1)^2 \> (u_0^2-c_1u_1^2)^N \> (du_0^2 -c_1du_1^2)
$$
$$
\nabla_{\bf x} = (N+1)^{-1} \> (u_0^2-c_1u_1^2)^{-N} \> \eta \>
                                 H({\bf u};c_1)^N \> \eta \> \nabla_{\bf u}
$$
$$
\partial_{x_0x_0} - c_1 
\partial_{x_1x_1} =
(N+1)^{-2} \> (u_0^2-c_1u_1^2)^{-N} \> 
(\partial_{u_0u_0}-c_1\partial_{u_1u_1})
\eqno (34)
$$
which hold outside the domain 
$\{(u_0,u_1) \in {\bf R}^2 \ | \ u_0^2 - c_1u_1^2=0\}$.

\sb

{\bf Example 2}. The mapping $T[-1;c_1,c_2; {\bf \varepsilon}]$, 
where ${\bf \varepsilon} = {\rm diag} (1,-1,1,1)$, 
yields the generic nonquadratic formulae
$$
x_0 = {1 \over \rho^2} 
\big ( u_0^2 + c_1 u_1^2 - c_2 u_2^2 + c_1c_2 u_3^2 \big )
$$
$$
x_1 = -{2 \over \rho^2}    u_0u_1 \qquad 
x_2 = -{2 \over \rho^2} c_1u_3u_1 \qquad 
x_3 = -{2 \over \rho^2}    u_2u_1 
\eqno (35)
$$
with $\rho^2 = u_0^2 - c_1 u_1^2 - c_2 u_2^2 + c_1c_2 u_3^2$. The general 
property (33) can be particularized as
$$
x_0^2 - c_1x_1^2 - c_2x_2^2 + c_1c_2x_3^2 = 1 
\eqno (36)
$$
Therefore, the nonquadratic transformation of type $B_{-1}$ with 
$2m = 4$ and $c_1=c_2=-1$
corresponds to the ${\bf R}^4 \rightarrow S^3$ Fock 
stereographic projection well known in quantum mechanics.

  
\vfill\eject

\noindent {\bf APPLICATIONS}

\sb

The applications of the Hurwitz transformations $T[N;c;{\bf \varepsilon}]$ range
from number theory (Hurwitz problem, Pythagorean and Diophantine equations)
to theoretical 
physics (classical and quantum mechanics, gauge theory). We give
below a few simple examples of application of some Hurwitz transformations. 

\sb

\noindent
{\bf Number Theory}

\sb

First, the Hurwitz theorem shows that Eq.~(6) admits solutions only for 
$n = 2$, $4$
and 8. (The case $n=1$ is trivial.) The noncompact extension of this theorem
concerns the expressions
$$
       \sum_{\alpha = 0}^{ m-1} w_{\alpha}^2 -
       \sum_{\alpha = m}^{2m-1} w_{\alpha}^2 =
\bigg( \sum_{\alpha = 0}^{ m-1} u_{\alpha}^2 - 
       \sum_{\alpha = m}^{2m-1} u_{\alpha}^2 \bigg) 
\bigg( \sum_{\alpha = 0}^{ m-1} v_{\alpha}^2 -
       \sum_{\alpha = m}^{2m-1} v_{\alpha}^2 \bigg) 
\eqno (37)
$$
for $2m = 2$, 4 and 8. As an interesting result, the solutions 
$w_\alpha = \phi_\alpha (u_\beta , v_\beta)$ of Eq.~(37) are given by 
${\bf w} = H({\bf u};c) {\bf v}$.

Second, the three generic Diophantine equations
$$
A^2-c_1B^2-c_2C^2+c_1c_2D^2-c_3E^2+c_1c_3F^2+c_2c_3G^2-c_1c_2c_3H^2 = I^2 
\eqno (38)
$$
$$
A^2 - c_2B^2 + c_1c_2C^2 - c_3D^2 + c_1c_3E^2 = F^2 
\eqno (39)
$$
and
$$
A^2-c_2B^2+c_1c_2C^2-c_3D^2+c_1c_3E^2+c_2c_3F^2-c_1c_2c_3G^2 = H^2 
\eqno (40)
$$
admit solutions 
$(A, B, \cdots, I) \in {\bf Z}^9$,
$(A, B, \cdots, F) \in {\bf Z}^6$ and 
$(A, B, \cdots, H) \in {\bf Z}^8$
corresponding to quadratic transformations of type
$A_1$, $B_1$ and $C_1$, respectively (see Table~2). 
In Eqs.~(38), (39) and (40), 
each $c_i$ ($i=1,2,3$) can take the values  $+1$
or $-1$ (as in what precedes) in the case $2m=8$. Then, the solutions are 
quadratic functions of $u_{\alpha} \in {\bf Z}$ ($\alpha =0,1,\cdots,7$) 
given by
Eqs.~(17) and (18) for (38),
Eqs.~(21) and (22) for (39) and 
Eqs.~(26) and (27) for (40). Other Diophantine equations are
obtained in the cases $2m=4$ and $2m=2$ by taking $c_3=0$ and $c_2=c_3=0$,
respectively. 

We have reported in Table 2 some solutions for the Diophantine
equations (38), (39) et (40) and their particular cases. In Table~2, the 
             solutions are $(A,B,\cdots,I) \in {\bf Z}^9$,
                           $(A,B,\cdots,F) \in {\bf Z}^6$ and 
                           $(A,B,\cdots,H) \in {\bf Z}^8$ for 
$c_i=\pm1$ with $i=1,2,3$. The two particular
Diophantine equations inside round brackets admit solutions 
$(A,B,C,D,I) \in {\bf Z}^5$ and $(A,B,C,F) \in {\bf Z}^4$ corresponding to 
$u_4=u_5=u_6=u_7$ and $c_3=0$ while the Diophantine equation in square brackets
admits solutions $(A,B,I) \in {\bf Z}^3$ 
corresponding to $u_2=u_3=u_4=u_5=u_6=u_7$ and $c_2=c_3=0$.

For example, for $c_1 = c_2 = c_3 - 1 = -1$ we get from Eq.~(39)
the Pythagorean equation (see also Table~2)
$$
A^2+B^2+C^2=F^2
\eqno (41)
$$
with the solutions
$$
A = u_0^2 + u_1^2 - u_2^2 - u_3^2 \qquad
B = 2(u_0u_2 - u_1u_3) 
$$
$$ 
C = 2(u_0u_3 + u_1u_2)            \qquad 
F = u_0^2 + u_1^2 + u_2^2 + u_3^2 
\eqno (42)
$$
Similarly, for $c_2=c_3=0$ we get from Eq.~(38) the Diophantine equation
$$
A^2-c_1B^2=I^2
\eqno (43)
$$
with the solutions
$$
A = u_0^2 + c_1u_1^2 \qquad 
B = 2u_0u_1          \qquad 
I = u_0^2 - c_1u_1^2 
\eqno (44)
$$
Note that the introduction of specific relations between the $u_\alpha$ in 
Eqs.~(17), (21) and (26) may lead to new Diophantine equations. 
For instance, by putting $u_1=0$ and $u_2=u_3$ in Eq.~(42), we obtain that
$$
A = u_0^2 - 2u_2^2 \qquad 
B = 2u_0u_2        \qquad 
F = u_0^2 + 2u_2^2 
\eqno (45)
$$
are solutions of
$$
A^2 + 2 B^2 = F^2
\eqno (46)
$$
an equation which comes from Eq.~(41). 

\sb

\centerline {
          {\bf Table 2.}~~Diophantine 
          equations in ${\bf Z}^r$ for $r=3,4,5,6,8,9$ 
          with solutions.} 

\def\init{\tabskip 0pt\offinterlineskip}
\def\crr{\cr\noalign{\hrule}}

$$\vbox{\init\halign to 14.2 true cm{\strut#&\vrule#%
\tabskip=0em plus 2 em&
\hfil$#$\hfil&
\vrule#&
\hfil$#$\hfil&
\vrule#\tabskip 0pt\crr
&& && &\cr
&&\hbox {Equation} &&\hbox {Solution} &\cr
&& && &\crr
&& && &\cr
&&\phantom{(}A^2-c_1B^2-c_2C^2\hfill 
&&A=u^2_0+c_1u^2_1+c_2u^2_2-c_1c_2u^2_3+c_3u^2_4\hfill&\cr

&&\phantom{(A^2}+c_1c_2D^2-c_3E^2+c_1c_3F^2\hfill 
&&\phantom{A=u^2_0}-c_1c_3u^2_5-c_2c_3u^2_6+c_1c_2c_3u^2_7\hfill &\cr

&&\phantom{(A^2}+c_2c_3G^2-c_1c_2c_3H^2=I^2\hfill 
&&B=2u_0u_1 \quad C=2u_0u_2 \quad D=2u_0u_3\hfill &\cr

&&(A^2-c_1B^2-c_2C^2\hfill 
&&E=2u_0u_4 \; F=2u_0u_5 \; G=2u_0u_6 \; H=2u_0u_7\hfill &\cr

&&\phantom{(A^2}+c_1c_2D^2=I^2)\hfill 
&&I=u_0^2-c_1u_1^2-c_2u_2^2+c_1c_2u^2_3\hfill &\cr

&&[A^2-c_1B^2=I^2]\hfill 
&&\phantom{I=u^2_6}-c_3u^2_4+c_1c_3u^2_5+c_2c_3u^2_6-c_1c_2c_3u^2_7\hfill &\cr
&& && &\crr
&& && &\cr

&& &&A=u^2_0-c_1u^2_1+c_2u^2_2-c_1c_2u^2_3+c_3u^2_4\hfill &\cr
&& &&\phantom{A=u^2_0}-c_1c_3u^2_5+c_2c_3u^2_6-c_1c_2c_3u^2_7\hfill &\cr

&&\phantom{(}A^2 - c_2B^2 + c_1c_2C^2\hfill 
&&B=2(u_0u_2+c_1u_1u_3+c_3u_4u_6-c_1c_3u_5u_7)\hfill &\cr

&&\phantom{(A^2}-c_3D^2+c_1c_3E^2=F^2\hfill 
&&C=2(u_0u_3+u_1u_2-c_3u_5u_6+c_3u_4u_7)\hfill &\cr

&&(A^2-c_2B^2+c_1c_2C^2=F^2)\hfill 
&&D=2(u_0u_4+c_1u_1u_5-c_2u_2u_6+c_1c_2u_3u_7)\hfill &\cr

&& &&E=2(u_0u_5+u_1u_4+c_2u_3u_6-c_2u_2u_7)\hfill &\cr
&& &&F=u_0^2-c_1u^2_1-c_2u^2_2+c_1c_2u^2_3\hfill &\cr
&& &&\phantom{F=u^2_0}-c_3u^2_4+c_1c_3u^2_5+c_2c_3u^2_6-c_1c_2c_3u^2_7\hfill 
&\cr

&& && &\crr
&& && &\cr

&& &&A=u^2_0-c_1u^2_1-c_2u^2_2+c_1c_2u^2_3-c_3u^2_4\hfill &\cr
&& &&\phantom{A=u^2_0}+c_1c_3u^2_5-c_2c_3u^2_6+c_1c_2c_3u^2_7\hfill &\cr
&& &&B=2(-c_3u_4u_6+c_1c_3u_5u_7)\hfill &\cr

&&A^2-c_2B^2+c_1c_2C^2\hfill 
&&C=2(-c_3u_4u_7+c_3u_5u_6)\hfill &\cr

&&\phantom{A^2} - c_3D^2 + c_1c_3E^2 + c_2c_3F^2\hfill 
&&D=2(c_2u_2u_6-c_1c_3u_3u_7)\hfill &\cr

&&\phantom{A^2}-c_1c_2c_3G^3=H^2\hfill 
&&E=2(c_2u_2u_7-c_2u_3u_6)\hfill &\cr

&& &&F=2(u_0u_6-c_1u_1u_7)\hfill &\cr

&& &&G=2(u_0u_7-u_1u_6)\hfill &\cr

&& &&H=u^2_0-c_1u^2_1-c_2u^2_2+c_1c_2u^2_3\hfill &\cr

&& &&\phantom{H=u^2_0}-c_3u^2_4+c_1c_3u^2_5+c_2c_3u^2_6-c_1c_2c_3u^2_7\hfill &\cr
&& && &\crr
}}$$

\sa 

Third, the Hurwitz transformations of type $A_N$, $B_N$ and $C_N$ with 
$N = 2,3,4,\cdots$ lead to Diophantine equations where the variables 
$A,B,C,\cdots$ 
are homogeneous polynomials of degree $N+1$ in $u_\alpha$ 
($\alpha = 0,1,\cdots,2m-1$). As
a typical example, from the Hurwitz transformation $T[N;c_1;{\bf 1}]$ 
we see that
$$
A^2 - c_1B^2 = C^{N+1} \qquad N = 2,3,4, \cdots 
\eqno (47)
$$
admit as (particular) solution $A$, $B$ and $C$ such that
$$  
  \pmatrix{ A\cr 
            B\cr}  =
  \pmatrix{ u_0 & c_1 u_1\cr
            u_1 &     u_0\cr}^N 
  \pmatrix{ u_0\cr
            u_1\cr}
\qquad
  C = (u_0 \ u_1) 
  \pmatrix{ 1 &    0\cr
            0 & -c_1\cr}
  \pmatrix{ u_0\cr
            u_1\cr}
\eqno (48)
$$
Particular solutions of more complicated equations generalizing Eq.~(47) 
can be obtained from transformations of type $A_N$, $B_N$ and $C_N$.

\sb

\noindent
{\bf Dynamical Systems} 

\sb

As a second illustration, we end up with a toy 
example concerning applications of Hurwitz transformations to dynamical 
systems. Let us consider the two-dimensional Schrvdinger equation
$$
\bigg [
- {1 \over 2}      \big (
 \partial_{x_0x_0} + 
 \partial_{x_1x_1} \big ) 
- {Z \over r^\alpha} \bigg ] \psi = E \psi \qquad
\alpha \in {\bf R}
\eqno (49)
$$
where $Z \in {\bf R}$ and $r= \sqrt {x_0^2 + x_1^2}$. 
The application to Eq.~(49) of a
transformation $T[N;-1;{\bf 1}]$ (of type $A_N$) with $N \in {\bf Z}-\{-1\}$
leads to the partial differential equation
$$
 \bigg [ - {1 \over 2} 
 \big ( \partial_{u_0u_0} + \partial_{u_1u_1} \big )
- \big ( N+1 \big )^2 E \rho^{2N} \bigg ]     {\widehat \psi}
= \big ( N+1 \big )^2 Z \rho^{2N-\alpha(N+1)} {\widehat \psi}
\eqno (50)
$$
where $\rho =\sqrt {u_0^2+u_1^2}$ and 
${\widehat \psi} \equiv {\widehat \psi}({\bf u})$ is the transform of 
$          \psi  \equiv           \psi ({\bf x})$ under the considered 
transformation. Furthermore, let us impose that $2N - \alpha(N+1) = 0$.
Then, Eq.~(50) reduces to 
$$
\bigg [ - {1 \over 2} \big ( \partial_{u_0u_0} + \partial_{u_1u_1} \big ) -
\big ( N+1 \big )^2 E \rho^{2N} \bigg ] {\widehat \psi} =
\big ( N+1 \big )^2 Z                   {\widehat \psi}
\eqno (51)
$$
so that the transformation of type $A_N$ allows to transform the 
${\bf R}^2$ Schrvdinger equation for the potential 
$- Z \big ( x_0^2 + x_1^2 \big )^{-N/(N+1)}$
and the energy $E$ into the ${\bf R}^2$ Schrvdinger equation for the 
potential
$- \big ( N+1 \big )^2 E \big (u_0^2 + u_1^2 \big )^N$
and the energy $(N+1)^2Z$. 
(Note that in such a transformation the rtles of the energy and the coupling
constant are interchanged.) The solutions for $\alpha \in {\bf Z}$ and 
$N \in {\bf Z} - \{ -1 \}$ correspond to the couples $(\alpha, N) = (1,1)$,
$(3,-3)$ and $(4,-2)$. (The solution $(0,0)$ is trivial~!) In other words, the
${\bf R}^2$ Schrvdinger equations for the potential $1/r$ (Coulomb potential), 
$1/r^3$ and
$1/r^4$ are transformed into the ${\bf R}^2$ 
Schrvdinger equations for the potentials
$\rho^2$ (harmonic oscillator potential), $1/\rho^6$ and $1/\rho^4$,
respectively. 

To close this paper, let us mention that there exist other applications to
dynamical systems. In classical mechanics, the Kustaanheimo-Stiefel
transformation (a transformation of type $B_1$) was used for the regularization
of the Kepler problem.$^3$ 
In quantum mechanics, the latter transformation made it
possible to transform the Schrvdinger equation for the three-dimensional
hydrogen atom (in an electromagnetic field) into a Schrvdinger equation for a
four-dimensional isotropic harmonic oscillator (with quartic and sextic
anharmonic terms) subject to a constraint.$^{4,5,8}$

\sb

\noindent
{\bf Acknowledgments}

\sb 

Part of the present work was done in collaboration with several 
people. In this respect, the author is grateful to his 
co-authors of Refs.~5 and 8-13 for fruitful collaborations. 
The applications to Number Theory were not published before. 
In this connection, the author acknowledges Profs. J.-P. 
Nicolas, T. Ono, and A. Schinzel for interesting comments. 
A preliminary version of this paper appears in Ref.~16. 
Finally, the author is very indebted to Prof.~Bruno Gruber 
to have given him the opportunity to present this work 
to the beautiful symposium ``Symmetries in Science IX''.

\sb
\sd

\noindent {\bf REFERENCES}

\sb

\item{1.} 
 Y. Choquet-Bruhat, C. DeWitt-Morette and M. Dillard-Bleick, 
 {\it Analysis, Manifolds and Physics} (North-Holland, Amsterdam, 1982). 

\item{2.}  
 Y. Choquet-Bruhat and C. DeWitt-Morette, {\it Analysis,
 Manifolds and Physics. Part II: 92 Applications} (North-Holland, Amsterdam,
 1989). 

\item{3.}  
 P. Kustaanheimo and E. Stiefel, {\it J. reine angew. Math.} 
 {\bf 218} (1965) 204. 

\item{4.}  
 M. Boiteux, {\it J. Math. Phys.} {\bf 23} (1982) 1311. 

\item{5.}  
 M. Kibler and T. N\'egadi, 
 {\it Lett. Nuovo Cimento} {\bf 39} (1984) 319.

\item{6.}  
 I. V. Polubarinov, {\it On Application of Hopf Fiber 
 Bundles in Quantum Theory}, Report E2-84-607 (JINR, Dubna, 1984). 

\item{7.}  
 T. Iwai, {\it J. Math. Phys.} {\bf 26} (1985) 885. 

\item{8.}  
 A. C. Chen and M. Kibler, {\it Phys. Rev.} 
 {\bf A31} (1985) 3960.

\item{9.}  
 M. Kibler, A. Ronveaux and T. N\'egadi, {\it J. Math. Phys.} 
 {\bf 27} (1986) 1541. 

\item{10.}  
 D. Lambert and M. Kibler, {\it J. Phys. A: Math. Gen.} 
 {\bf 21 } (1988) 307. 

\item{11.}  
 M. Kibler and P. Winternitz, {\it J. Phys. A: Math. Gen.} 
 {\bf 21 } (1988) 1787. 

\item{12.}  
 M. Kibler and P. Labastie, 
 in {\it Group Theoretical Methods in  Physics}, 
 eds. Y. Saint-Aubin and L. Vinet (World Scientific, Singapore, 1989).

\item{13.}  
 M. Hage Hassan and M. Kibler, in 
 {\it Le probl\`eme de factorisation de Hurwitz}, 
 eds. A. Ronveaux and D. Lambert (FUNDP, Namur, 1991). 

\item{14.}  
 L. E. Dickson, {\it Ann. Math.} {\bf 20} (1919) 155. 

\item{15.}  
 T. Ono, {\it Variations on a Theme of Euler: Quadratic
 Forms, Elliptic Curves, and Hopf Maps} (Plenum Press, New York, 1994).

\item{16.}  
 M. Kibler, in 
 {\it Proceedings of the 
 IIIrd Catalan Days on Applied Mathematics}, 
 eds. J. Chavarriga and J. Gin\'e (Institut d'Estudis Ilerdencs, Lleida, 1996).

\bye